\documentclass[11pt,letterpaper]{article}
\usepackage[margin=1in]{geometry}
\usepackage[T1]{fontenc}
\usepackage{lmodern}
\usepackage{amsmath,amssymb,amsthm}
\usepackage{microtype}
\usepackage{booktabs}
\usepackage{tikz}
\usepackage{needspace}
\usepackage[hyphens]{url}
\usepackage[hidelinks]{hyperref}
\usepackage{doi}

\theoremstyle{plain}
\newtheorem{theorem}{Theorem}
\newtheorem{lemma}[theorem]{Lemma}
\newtheorem{corollary}[theorem]{Corollary}
\theoremstyle{definition}
\newtheorem{definition}[theorem]{Definition}

\newtheorem{example}[theorem]{Example}
\newcommand{\floor}[1]{\left\lfloor #1\right\rfloor}
\newcommand{\ceil}[1]{\left\lceil #1\right\rceil}
\newcommand{\coef}[1]{[z^{#1}]}
\hypersetup{
  pdftitle={Halving Partitions and Hyperbinary Expansions},
  pdfauthor={Max A. Alekseyev},
  pdfsubject={Halving partitions, hyperbinary expansions, and Stern's sequence}
}

\title{Halving Partitions and Hyperbinary Expansions}
\author{Max A. Alekseyev\thanks{The George Washington University,
  Washington, DC, United States.
  Email: \href{mailto:maxal@gwu.edu}{\nolinkurl{maxal@gwu.edu}}.}}
\date{}

\begin{document}
\pagenumbering{arabic}
\maketitle
\begin{abstract}
Starting from a positive integer, repeatedly halve the current value,
rounding up or down.  Which resulting sequences have terms that add to
the starting integer?  A simple balance identity turns the rounding
choices into hyperbinary expansions, in which each power of two occurs
at most twice.  This correspondence gives exact counts by length and an
injection into the expansions counted by Stern's diatomic sequence.  We
characterize equality in the resulting bound by the binary digits of the
odd part of the integer.  We also determine the unique longest partition
and give a short recursive counting algorithm.  When the required total
is allowed to vary, the chain counts eventually equal Stern's number,
and we determine the exact threshold.  A direct proof shows that the
counts of hyperbinary expansions in the three length classes modulo
three differ by at most one.  For sufficiently long halving partitions,
the counts by length first weakly increase and then weakly decrease,
although the full sequence need not do so.
\end{abstract}
\smallskip

\section{Halving partitions}\label{sec:intro}

For a positive integer $m$, write
\[
  \mathcal H(m):=\left\{\floor{m/2},\ceil{m/2}\right\}.
\]
When $m$ is even, the two displayed values coincide, so $\mathcal H(m)$ is a
set with one element.  Clark Kimberling introduced the following notion in the OEIS entry
A349552 \cite{OEIS}.

\begin{definition}
A \emph{halving partition} of a positive integer $n$ is a partition
\[
  n=p_1+p_2+\cdots+p_k
\]
with positive parts such that
\[
  p_1\in\mathcal H(n),\qquad
  p_i\in\mathcal H(p_{i-1})\quad(2\le i\le k).
\]
We also regard the empty partition as the unique halving partition of $0$.
For $n,k\ge0$, let $a_k(n)$ denote the number of halving partitions of $n$
with exactly $k$ parts.  Thus $a_0(0)=1$ and $a_k(0)=0$ for $k>0$.
For every $n\ge0$ we put
\[
  a(n):=\sum_{k\ge0}a_k(n).
\]
\end{definition}

For $n\ge1$, let
$L(n):=\max\{k\ge0:a_k(n)>0\}$ denote the maximum number of parts in a
halving partition of $n$.

The first values of $a(n)$ for $n\ge0$, sequence A349552, are
\[
1,1,1,2,1,3,2,3,1,3,3,5,2,5,3,4,1,4,3,6,\ldots.
\]
For example,
\[
 13=7+4+2=6+3+1+1+1+1
\]
gives halving partitions of two different lengths.  All five halving
partitions of $13$ appear in Table~\ref{tab:13} below.
Each part is constrained by the preceding part, but the sum condition
$p_1+\cdots+p_k=n$ involves all the parts.  The counting problem is to
satisfy both requirements at once.

Two elementary observations already reveal the role of binary digits.
\begin{lemma}\label{prop:dyadic}
For $n\ge1$ and $k\ge2$,
\begin{equation}\label{eq:dyadic-length}
 a_k(2n)=a_{k-1}(n).
\end{equation}
Moreover $a_1(2n)=0$.  In particular,
\begin{equation}\label{eq:dyadic}
 a(2n)=a(n).
\end{equation}
\end{lemma}

\begin{proof}
Every halving partition of $2n$ must begin with $n$, since floor and ceiling
coincide.  Removing that first part leaves a halving partition of $n$, because
the remaining parts sum to $n$ and obey the same successive-halving rule.
Conversely, prepend $n$ to any halving partition of $n$.
\end{proof}

Thus $a(n)$ depends only on the odd part of $n$, that is, its largest odd
divisor.  In particular, $a(2^m)=1$, as recorded in A349552.  We next
consider the maximum number of parts.  Write $s_2(n)$ for the number of
ones in the ordinary binary expansion of $n$.  Throughout, binary
expansions of positive integers are written without leading zeros.
\begin{theorem}\label{thm:longest}
For $n\ge1$,
\begin{equation}\label{eq:maxlength}
  L(n)=\floor{\log_2 n}+s_2(n),
\end{equation}
and exactly one halving partition has $L(n)$ parts.
\end{theorem}
\begin{proof}
The case $n=1$ is immediate, so assume $n\ge2$ and put $h:=\floor{\log_2 n}$.  If a halving partition has
$k<h$ parts, the claimed upper bound is automatic.  Suppose therefore that
$k\ge h$.  For $1\le i\le h$, every possible $p_i$ satisfies
$p_i\ge\floor{n/2^i}$.  The classical binary identity
\begin{equation}\label{eq:floor-sum}
 \sum_{i=1}^{h}\floor{\frac{n}{2^i}}=n-s_2(n)
\end{equation}
shows that the first $h$ parts have sum at least $n-s_2(n)$.
Indeed, each nonzero binary digit in position $j$ contributes
$\sum_{i=1}^{j}2^{j-i}=2^j-1$ to the left side of \eqref{eq:floor-sum},
with an empty sum when $j=0$.  Any further parts are positive, so there can be at
most $s_2(n)$ of them.  Hence $k\le h+s_2(n)$.

Equality forces equality at every step in the preceding estimates: the first
$h$ parts must be the floors displayed above, and every remaining part must be
$1$.  Their sum is $n$ by \eqref{eq:floor-sum}.  The displayed floors
are obtained by rounding down at each step, and once $1$ is reached it
may be repeated by rounding up.  Thus the candidate is indeed a halving partition, and it
is unique.
\end{proof}

The total count is therefore the finite sum
\begin{equation}\label{eq:finite-total}
 a(n)=\sum_{k=0}^{L(n)}a_k(n).
\end{equation}
Restrictions involving halving or doubling also occur in Cayley
compositions and non-squashing partitions
\cite{KonvalinkaPak,SloaneSellers}.  Another comparison comes from
payphone permutations, which describe successive choices of vacant phones
that maximize the distance to the nearest occupied phone \cite{AlekseyevPayphone}.
Subtracting one from successive halving values and reversing their order
gives the moves $y\mapsto2y,2y+1,2y+2$.  The payphone interval calculation
uses the same moves, with the middle move counted twice.  Halving
partitions also require the chosen parts to sum to $n$;
Appendix~\ref{sec:appendix} gives the detailed comparison.
Liu and Wang have also studied payphone permutations through a different
construction: rooted trees that encode permutations
\cite{LiuWangPermutationTrees}.

Section~\ref{sec:correspondence} encodes the rounding choices by digits
and turns the prescribed sum into a balance identity.
Theorem~\ref{thm:main} uses this identity to count partitions of each
length by hyperbinary expansions, in which each power of two occurs at
most twice.  Theorem~\ref{prop:signed} gives a signed-digit description
and the largest integer admitting a partition of any prescribed
positive length.
Section~\ref{sec:stern} then identifies which expansions arise from
halving partitions.  Theorem~\ref{thm:deficit} gives an injection into
the hyperbinary expansions of $n-1$, bounding $a(n)$ by the corresponding
term of Stern's diatomic sequence.  Corollary~\ref{cor:equality}
characterizes equality by the presence of at most one zero in the binary
expansion of the odd part of $n$.

Section~\ref{sec:divide} allows the total of a halving chain to vary.
Theorem~\ref{thm:Arec} gives a recurrence for these auxiliary counts,
and Theorem~\ref{thm:stabilization} determines the exact threshold beyond
which they equal Stern's term.
Corollary~\ref{cor:binary-stabilization} shows that exactly two binary
zeros in the odd part of $n$ make $a(n)$ one less than Stern's term.
Section~\ref{sec:fences} studies the distribution by length:
Theorem~\ref{prop:trisection} shows that hyperbinary expansions of any
fixed integer are distributed among the length classes modulo three
with counts differing by at most one.  A correspondence with certain
partially ordered sets, called fences, shows in
Corollary~\ref{cor:unimodal} that the counts for sufficiently long
partitions are unimodal: they first weakly increase and then weakly
decrease.  Example~\ref{ex:nonunimodal} shows why the sequence of counts
for all lengths need not be unimodal.  Section~\ref{sec:discussion} discusses further questions,
and Appendix~\ref{sec:appendix} collects supplementary formulas and the
detailed payphone comparison.

\section{The balance identity and hyperbinary expansions}\label{sec:correspondence}
Set $p_0:=n$.  Encode each halving step by the unique digit
$q_i\in\{0,1,2\}$ satisfying
\begin{equation}\label{eq:qdef}
 p_{i-1}=2p_i+q_i-1.
\end{equation}
Thus $q_i=1$ when $p_{i-1}$ is even.  When $p_{i-1}$ is odd, rounding
up gives $q_i=0$, and rounding down gives $q_i=2$.  In these digits, the
requirement that the parts sum to $n$ is especially simple.
\Needspace{12\baselineskip}
\begin{lemma}\label{lem:balance}
For a halving partition $p_1+\cdots+p_k=n$ with $n\ge1$,
\begin{equation}\label{eq:balance}
 \sum_{i=1}^kq_i=k-p_k,
\end{equation}
and
\begin{equation}\label{eq:qbinary}
 n-1=(p_k-1)2^k+\sum_{i=1}^kq_i2^{i-1}.
\end{equation}
\end{lemma}
\begin{proof}
From \eqref{eq:qdef},
$p_{i-1}-p_i=p_i+q_i-1$.  Summing over $i$ gives
\[
 n-p_k=\sum_{i=1}^k p_i+\sum_{i=1}^k(q_i-1)
       =n+\sum_{i=1}^kq_i-k,
\]
and therefore $\sum_iq_i=k-p_k$, which is \eqref{eq:balance}.

Finally, iterating $p_{i-1}=2p_i+q_i-1$ yields
\[
 n=2^kp_k+\sum_{i=1}^k(q_i-1)2^{i-1}
   =2^kp_k+\sum_{i=1}^kq_i2^{i-1}-(2^k-1),
\]
which gives \eqref{eq:qbinary}.
\end{proof}

The powers of two in \eqref{eq:qbinary} lead to a classical digit system.
A \emph{hyperbinary expansion} of a nonnegative integer $m$ is an
expression
\[
 m=\sum_{j\ge0}e_j2^j,\qquad e_j\in\{0,1,2\};
\]
each power of two is used at most twice
\cite{Reznick,DilcherEricksen}.  Its \emph{length} is the number of
summands, $\sum_j e_j$, not the number of digit positions.  Define the
hyperbinary length polynomial
\begin{equation}\label{eq:Hdef}
  H_m(z):=
  \sum_{\substack{e_j\in\{0,1,2\}\\
                   \sum_{j\ge0}e_j2^j=m}}
  z^{\sum_{j\ge0} e_j}.
\end{equation}
For a formal power series $F(z)$, we write $\coef{r}F(z)$ for the coefficient
of $z^r$ in $F(z)$, with $\coef{r}F(z):=0$ for $r<0$.  Thus
$\coef{r}H_m(z)$ counts hyperbinary expansions of $m$ having $r$ summands.
The length polynomials are collected in the two-variable generating function
\begin{equation}\label{eq:bivariate-H}
  G(x,z):=\sum_{m\ge0}H_m(z)x^m
  =\prod_{j\ge0}\left(1+zx^{2^j}+z^2x^{2^{j+1}}\right).
\end{equation}

Factoring out the $j=0$ term gives
\begin{equation}\label{eq:Gfunc}
 G(x,z)=(1+zx+z^2x^2)G(x^2,z).
\end{equation}
Comparing odd and even powers of $x$ yields
\begin{equation}\label{eq:Hrec}
 H_{2m+1}(z)=zH_m(z),\qquad
 H_{2m}(z)=H_m(z)+z^2H_{m-1}(z),
\end{equation}
for $m\ge0$, with $H_0(z)=1$ and $H_{-1}(z):=0$.
\Needspace{7\baselineskip}
Putting $z=1$ counts the expansions without distinguishing their lengths.  Their
classical generating function is
\begin{equation}\label{eq:stern-product}
 \prod_{j\ge0}(1+x^{2^j}+x^{2^{j+1}})
 =\sum_{m\ge0}s(m+1)x^m,
\end{equation}
where Stern's diatomic sequence is defined by
\[
\begin{gathered}
 s(0):=0,\qquad s(1):=1,\\
 s(2m):=s(m),\qquad s(2m+1):=s(m)+s(m+1)\quad(m\ge1).
\end{gathered}
\]
In particular,
\begin{equation}\label{eq:sternH}
 H_m(1)=s(m+1).
\end{equation}
The length polynomial also plays a central role in the recent work of
McConville, Propp, and Sagan \cite{McConvilleProppSagan}.

\subsection{From partitions to expansions}
Let $P=(p_1,\ldots,p_k)$ be a halving partition of $n\ge1$.  Write
$p_k-1=\sum_{j\ge0}b_j2^j$ in ordinary binary.  Then
\eqref{eq:qbinary} becomes
\[
 n-1=\sum_{i=1}^kq_i2^{i-1}+\sum_{j\ge0}b_j2^{k+j}.
\]
The digits in positions $0,\ldots,k-1$ are $q_1,\ldots,q_k$; those in
positions $k,k+1,\ldots$ are the binary digits of $p_k-1$.
Together they form a hyperbinary expansion of $n-1$, denoted by $\Phi_n(P)$.
For $13=7+4+2$, the first three digits are $(0,0,1)$ and $p_3-1=1$, giving
$12=8+4$.  Table~\ref{tab:13} shows all five correspondences.

\begin{table}[ht]
\centering
\caption{The five halving partitions of $13$ and their images under $\Phi_{13}$.}
\label{tab:13}
\begin{tabular}{cc}
\toprule
Halving partition of $13$ & Hyperbinary expansion of $12$\\
\midrule
$7+4+2$ & $8+4$\\
$7+3+2+1$ & $8+2+2$\\
$7+3+1+1+1$ & $4+4+2+2$\\
$6+3+2+1+1$ & $8+2+1+1$\\
$6+3+1+1+1+1$ & $4+4+2+1+1$\\
\bottomrule
\end{tabular}
\end{table}

We will prove that $\Phi_n$ is injective and determine its image in
Section~\ref{sec:stern}.  First we adapt it to count partitions of a fixed
length.

\subsection{Counting by length}
For a prescribed length $k$, we subtract a multiple of $2^k$ from
\eqref{eq:qbinary} so that the remaining expansion has a prescribed
length as well.  Define, for $n,k\ge1$,
\begin{equation}\label{eq:QMdef}
 Q_k(n):=\max\left\{1,\floor{\frac{n}{2^k}}\right\},
 \qquad
 M_k(n):=n-1-(Q_k(n)-1)2^k.
\end{equation}
Writing $h:=\floor{\log_2 n}$, the second quantity can also be written as
\begin{equation}\label{eq:Mcases}
 M_k(n)=
 \begin{cases}
   (n\bmod 2^k)+2^k-1, & k<h,\\
   n-1, & k\ge h.
 \end{cases}
\end{equation}
Indeed, for $k<h$ we have $Q_k(n)=\floor{n/2^k}$, so
$n-Q_k(n)2^k=n\bmod 2^k$; once $k\ge h$, we have $Q_k(n)=1$.
This gives the fixed-length correspondence.

\begin{theorem}\label{thm:main}
For $n,k\ge1$,
\begin{equation}\label{eq:main}
  a_k(n)=\coef{k-Q_k(n)}H_{M_k(n)}(z).
\end{equation}
Equivalently, $a_k(n)$ is the coefficient of
$x^{M_k(n)}z^{k-Q_k(n)}$ in $G(x,z)$.
\end{theorem}
\begin{lemma}\label{lem:lastpart}
After $k$ successive halvings starting from $n$, the last positive value belongs
to
\[
 \left\{\floor{n/2^k},\ceil{n/2^k}\right\},
\]
with $0$ excluded if it occurs.  Consequently, for every halving partition
of length $k$,
\begin{equation}\label{eq:lastpartQ}
 p_k-Q_k(n)\in\{0,1\}.
\end{equation}
\end{lemma}

\begin{proof}
First allow zero among the attainable values and argue by induction.  After $j$ halvings, every
attainable value belongs to
$\{\floor{n/2^j},\ceil{n/2^j}\}$; applying one further floor- or
ceiling-halving step again produces a value in
$\{\floor{n/2^{j+1}},\ceil{n/2^{j+1}}\}$.  If the lower value is $0$,
positivity leaves only $1$.  This is exactly \eqref{eq:lastpartQ} with
$Q_k(n)$ from \eqref{eq:QMdef}.
\end{proof}

We can now prove Theorem \ref{thm:main}.

\begin{proof}[Proof of Theorem \ref{thm:main}]
Let $p_1+\cdots+p_k=n$ be a halving partition.  Subtract
$(Q_k(n)-1)2^k$ from \eqref{eq:qbinary}.  We obtain
\begin{equation}\label{eq:Mexp}
 M_k(n)=(p_k-Q_k(n))2^k+\sum_{i=1}^kq_i2^{i-1}.
\end{equation}
By Lemma \ref{lem:lastpart}, the coefficient of $2^k$ is $0$ or $1$, while
each lower coefficient $q_i$ lies in $\{0,1,2\}$.  Thus
\eqref{eq:Mexp} is a hyperbinary expansion of $M_k(n)$.  Its number of
summands is, by the balance identity,
\[
 (p_k-Q_k(n))+\sum_{i=1}^kq_i
 =p_k-Q_k(n)+k-p_k=k-Q_k(n).
\]
So every halving partition counted by $a_k(n)$ produces one expansion counted
by the right side of \eqref{eq:main}.

Conversely, $M_k(n)<2^{k+1}$, so a hyperbinary expansion of $M_k(n)$ has no
nonzero digit above $2^k$, and its $2^k$-digit is necessarily $0$ or $1$.
Write that digit as $e_k$, and write the lower digits as
$q_i=e_{i-1}$ for $1\le i\le k$.  Set $p_k:=Q_k(n)+e_k$ and define
backwards
\[
 p_{i-1}:=2p_i+q_i-1.
\]
The numerical identity \eqref{eq:Mexp} gives $p_0=n$.  Since $p_k\ge1$
and $p_{i-1}=2p_i+q_i-1\ge p_i\ge1$, all reconstructed parts are positive and
nonincreasing.  Moreover $q_i-1\in\{-1,0,1\}$, so each $p_i$ is obtained by rounding up or down after halving $p_{i-1}$.  The assumed number of summands is
$k-Q_k(n)$, hence
\[
 \sum_{i=1}^kq_i=k-p_k,
\]
and reversing the calculation in Lemma \ref{lem:balance} shows that
$\sum_i p_i=n$.  Thus we recover a halving partition, uniquely.
\end{proof}
The piecewise formula \eqref{eq:Mcases} distinguishes two ranges of
lengths.  For $n\ge2$, put
\begin{equation}\label{eq:hdef}
 h:=\floor{\log_2 n}.
\end{equation}
For $k\ge h$ we have $Q_k(n)=1$ and $M_k(n)=n-1$.  Hence
\begin{equation}\label{eq:longlengths}
 a_k(n)=\coef{k-1}H_{n-1}(z)\qquad(k\ge h).
\end{equation}
Thus, for $k\ge h$, the count is a coefficient of the same polynomial
$H_{n-1}(z)$, with the exponent one less than the partition length.
For $k<h$, the integer $M_k(n)=(n\bmod2^k)+2^k-1$ depends on $k$.
For these shorter lengths, \eqref{eq:rescaled-ak} in
Appendix~\ref{sec:appendix} gives an alternative coefficient formula
in which the residue $n\bmod2^k$ appears only in the exponent of $x$.

\begin{example}\label{ex:13}
Take
\[
  13=6+3+2+1+1,
\]
so $k=5$ and, with $p_0=13$, the digits determined by
\eqref{eq:qdef} are
\[
 (q_1,\ldots,q_5)=(2,1,0,1,0).
\]
Here $Q_5(13)=1$ and $M_5(13)=12$.  Formula \eqref{eq:Mexp} becomes
\[
 12=0\cdot2^5+0\cdot2^4+1\cdot2^3+0\cdot2^2
       +1\cdot2+2\cdot1
   =8+2+1+1.
\]
This hyperbinary expansion has four summands, exactly
$k-Q_5(13)=4$.

Conversely, start from the hyperbinary digits
$(e_0,\ldots,e_5)=(2,1,0,1,0,0)$.  They give $p_5=Q_5(13)+e_5=1$ and
$q_i=e_{i-1}$.  Successive applications of the backward recurrence give
\[
 (p_5,p_4,p_3,p_2,p_1,p_0)=(1,1,2,3,6,13),
\]
recovering the original partition.

For comparison,
\[
  H_{12}(z)=z^2+z^3+2z^4+z^5.
\]
Thus the five hyperbinary expansions of $12$ have lengths $2,3,4,4,5$,
and Theorem \ref{thm:main} turns them into the five halving partitions of
$13$ with lengths $3,4,5,5,6$ shown in Table~\ref{tab:13}.
\end{example}

\subsection{A signed-digit formula and the short lengths}
The same balance identity gives a description that does not use generating functions.
\begin{theorem}\label{prop:signed}
For $n,k\ge1$,
\begin{equation}\label{eq:signed-ak}
 a_k(n)=\#\left\{(c_0,\ldots,c_{k-1})\in\{-1,0,1\}^k:
 \ \sum_{i=0}^{k-1}c_i>0,\quad
 n=\sum_{i=0}^{k-1}c_i(2^k-2^i)\right\}.
\end{equation}
Consequently,
\begin{equation}\label{eq:fixedlength-support}
 \max\{n\ge1:a_k(n)>0\}=(k-1)2^k+1,
 \qquad
 a_k\bigl((k-1)2^k+1\bigr)=1.
\end{equation}
\end{theorem}

\begin{proof}
Given a halving partition of length $k$, set
$c_{i-1}:=1-q_i\in\{-1,0,1\}$.  The balance identity gives
\[
 p_k=k-\sum_{i=1}^kq_i=\sum_{i=0}^{k-1}c_i>0.
\]
Substituting this and $q_i-1=-c_{i-1}$ into the iterated identity in the
proof of Lemma~\ref{lem:balance} yields
\[
 n=2^k\sum_{i=0}^{k-1}c_i-\sum_{i=0}^{k-1}c_i2^i
   =\sum_{i=0}^{k-1}c_i(2^k-2^i).
\]
Conversely, start with a sequence of signed digits counted on the right side of
\eqref{eq:signed-ak}.  Put
$p_k:=\sum_i c_i\ge1$ and $q_i:=1-c_{i-1}$, and reconstruct backwards by
$p_{i-1}:=2p_i+q_i-1$.  Then $p_{i-1}=2p_i-c_{i-1}\ge p_i\ge1$, so all parts
are positive and each step halves the preceding part and rounds up or down.  The signed
identity gives $p_0=n$, while
$\sum_iq_i=k-p_k$; reversing the balance calculation therefore gives
$\sum_i p_i=n$.  This proves the bijection.

Finally, every weight $2^k-2^i$ in \eqref{eq:signed-ak} is positive and each
$c_i\le1$.  Hence
\[
 n\le\sum_{i=0}^{k-1}(2^k-2^i)
   =k2^k-(2^k-1)=(k-1)2^k+1.
\]
Equality forces $c_0=\cdots=c_{k-1}=1$, which is admissible and unique.
\end{proof}

The largest integer in \eqref{eq:fixedlength-support} corresponds to
rounding up at every step.  The bound will also show that only a few
lengths below $\floor{\log_2 n}$ can contribute to $a(n)$.

\subsection{A binary formula for the coefficients}
The coefficients of $H_m$ can themselves be described using ordinary
binary digits.  Call two nonnegative integers \emph{binary-disjoint} if
they have no position at which both binary digits are one.
\begin{lemma}\label{prop:Hexplicit}
For $m\ge0$,
\begin{equation}\label{eq:Hexplicit}
 H_m(z)=
 \sum_{\substack{0\le r\le m/2\\
                  r\text{ and }m-2r\text{ are binary-disjoint}}}
 z^{s_2(m-2r)+2s_2(r)}.
\end{equation}
\end{lemma}

\begin{proof}
Given a hyperbinary expansion $m=\sum_j e_j2^j$, collect the powers used twice
and those used once:
\[
 r:=\sum_{e_j=2}2^j,
 \qquad
 t:=\sum_{e_j=1}2^j.
\]
Then $m=2r+t$, and $r$ and $t$ are binary-disjoint.  Conversely, every decomposition
$m=2r+t$ with $r$ and $t$ binary-disjoint determines a unique hyperbinary expansion:
the positions of the ones in $r$ mark the powers used twice, and those in $t$ mark the
powers used once.  Its number of summands is
$2s_2(r)+s_2(t)=2s_2(r)+s_2(m-2r)$, proving \eqref{eq:Hexplicit}.
\end{proof}

Combining Lemma \ref{prop:Hexplicit} with Theorem \ref{thm:main} gives a
fully explicit description of A349552.

\begin{corollary}\label{cor:explicit}
For $n\ge1$, $a(n)$ is the number of pairs of integers $(k,r)$ satisfying
\begin{align*}
 1&\le k\le L(n),\\
 0&\le r\le M_k(n)/2,\\
 r\text{ and }M_k(n)-2r&\text{ are binary-disjoint},\\
 s_2(M_k(n)-2r)+2s_2(r)&=k-Q_k(n).
\end{align*}
Here $Q_k(n)$ and $M_k(n)$ are defined in \eqref{eq:QMdef}.
\end{corollary}

\begin{proof}
For each $k$, Theorem \ref{thm:main} asks for
$\coef{k-Q_k(n)}H_{M_k(n)}(z)$; Lemma \ref{prop:Hexplicit} describes
exactly which $r$ contribute to that coefficient.  The upper limit $L(n)$ is the maximum possible length by definition and
Theorem~\ref{thm:longest}.
\end{proof}

\section{Stern's bound and its equality cases}\label{sec:stern}
We now determine which hyperbinary expansions arise from $\Phi_n$.
For a partition of length $k$, positions $0,\ldots,k-1$ record the
rounding choices, and positions $k,k+1,\ldots$ record $p_k-1$ in binary.
For a hyperbinary expansion $e=(e_0,e_1,\ldots)$ of $N$, a
\emph{cut at position $k$} means the division between digits $e_{k-1}$ and
$e_k$.  For $k\ge1$, define
\begin{equation}\label{eq:cut}
\begin{gathered}
 L_k(e):=\sum_{j<k}e_j,\qquad
 R_k(e):=\sum_{j\ge k}e_j2^{j-k},\\
 F_e(k):=1+L_k(e)+R_k(e)-k.
\end{gathered}
\end{equation}
Here $L_k(e)$ is the sum of the digits below the cut, while $R_k(e)$
is the integer represented by the remaining digits after shifting them
$k$ positions.  Call a cut at position $k$ \emph{eligible} if no digit
$2$ occurs at or above position $k$, and \emph{admissible} if it is eligible
and $F_e(k)=0$.  The two conditions require, respectively, an ordinary
binary expansion of $p_k-1$ and the balance identity.
\begin{lemma}\label{prop:cut}
For $n\ge1$, a hyperbinary expansion $e$ of $n-1$ lies in the image of $\Phi_n$ if and
only if it has an admissible cut.  The cut, when it exists, is unique.  If its
position is $k$, the unique preimage under $\Phi_n$ is a halving partition of
length $k$.
\end{lemma}

\begin{proof}
Starting with a halving partition of length $k$, expand $p_k-1$ in ordinary
binary and place those digits at positions $k,k+1,\ldots$ in
\eqref{eq:qbinary}; place $q_i$ at position $i-1$.  This gives a hyperbinary
expansion $e$ of $n-1$ with no digit $2$ at or above the cut.  Moreover
$L_k=\sum_iq_i$, $R_k=p_k-1$, and the balance identity gives
\[
 F_e(k)=1+(k-p_k)+(p_k-1)-k=0.
\]

Conversely, suppose $e$ has an admissible cut at $k$.  Put
\[
 p_k:=R_k(e)+1,\qquad q_i:=e_{i-1}\quad(1\le i\le k),
\]
and reconstruct backwards by $p_{i-1}:=2p_i+q_i-1$.  Because the digits at and above the cut
are binary, $R_k(e)$ is the ordinary binary integer represented above the cut,
so the numerical identity gives $p_0=n$.  Also $p_i\ge1$ and
$p_{i-1}\ge p_i$, while $q_i-1\in\{-1,0,1\}$, so the successive parts obey the
halving rule.  Finally, $F_e(k)=0$ says
\[
 \sum_{i=1}^k q_i=L_k(e)=k-1-R_k(e)=k-p_k,
\]
and reversing the balance calculation gives $\sum_i p_i=n$.

For uniqueness, at every cut the identities
$R_k=e_k+2R_{k+1}$ and $L_{k+1}=L_k+e_k$ give
\begin{equation}\label{eq:Fdrop}
 F_e(k+1)-F_e(k)=-R_{k+1}(e)-1<0.
\end{equation}
Thus $F_e$ is strictly decreasing and can vanish at most once.
\end{proof}
To state the deficit formula, we first narrow the range of short lengths.
Theorem~\ref{prop:signed} gives the exact upper bound on $n$ for a fixed length:
$a_k(n)=0$ for $n>(k-1)2^k+1$.  For $n\ge2$, if
$k<h:=\floor{\log_2 n}$ and we write $k=h-j$, then
\[
 2^h\le n\le(h-j-1)2^{h-j}+1
\]
forces $2^j+j+1\le h$.  Define
\begin{equation}\label{eq:Jdef}
 J(h):=\max\bigl(\{0\}\cup\{j\ge1:2^j+j+1\le h\}\bigr).
\end{equation}
Thus among the lengths $k<h$, only
$k=h-1,h-2,\ldots,h-J(h)$ can contribute.  In particular, there are at most
$J(h)\le\floor{\log_2 h}=O(\log\log n)$ potentially nonzero short lengths.
\Needspace{15\baselineskip}
\begin{theorem}\label{thm:deficit}
For every $n\ge1$, the map $\Phi_n$ from halving partitions of $n$ to
hyperbinary expansions of $n-1$ is injective.  Consequently,
\[
 a(n)\le s(n).
\]
Define $d(n):=s(n)-a(n)$ for $n\ge1$.  It is exactly the number of
hyperbinary expansions of $n-1$ lying outside the image of $\Phi_n$.
For $n\ge2$, put $h:=\floor{\log_2 n}$.  Then
\begin{equation}\label{eq:deficit}
 d(n)=\coef{h-2}\frac{H_{n-1}(z)}{1-z}
      -\sum_{j=1}^{J(h)}a_{h-j}(n).
\end{equation}
Empty sums are understood to be zero.
\end{theorem}
\begin{proof}
Lemma \ref{prop:cut} shows that $\Phi_n$ is injective and that its
image consists precisely of the hyperbinary expansions of $n-1$ possessing an
admissible cut.  Since all hyperbinary expansions of $n-1$ are counted by
$s(n)$, we obtain $a(n)\le s(n)$, and $d(n)$ counts exactly the expansions
outside the image of $\Phi_n$.

For the coefficient formula, assume $n\ge2$ and put $h:=\floor{\log_2 n}$.  Equation
\eqref{eq:longlengths} gives
\[
 a(n)=\sum_{k=1}^{h-1}a_k(n)
       +\sum_{r\ge h-1}\coef{r}H_{n-1}(z),
\]
whereas
\[
 s(n)=H_{n-1}(1)=\sum_{r\ge0}\coef{r}H_{n-1}(z).
\]
Hence
\[
 d(n)=\sum_{r=0}^{h-2}\coef{r}H_{n-1}(z)
      -\sum_{k=1}^{h-1}a_k(n).
\]
The first sum is $\coef{h-2}H_{n-1}(z)/(1-z)$.  By the restriction on
short lengths preceding Theorem~\ref{thm:deficit}, the second sum is
$\sum_{j=1}^{J(h)}a_{h-j}(n)$.  This proves \eqref{eq:deficit}.
\end{proof}

Substituting \eqref{eq:main} for the short-length terms expresses the
deficit entirely in terms of hyperbinary coefficients; see
\eqref{eq:deficit-explicit} in Appendix~\ref{sec:appendix}.

For example, $H_8(z)=z+z^2+z^3+z^4$, whereas the halving partitions of
$9$ have lengths $3,4,5$.  The single-term expansion $8$ is omitted by
$\Phi_9$: its cut function takes the values $1$ and $-1$ at positions
$2$ and $3$, so it skips zero.  Hence $a(9)=3<s(9)=4$.

Equality in Stern's bound has a simple binary description.
\begin{corollary}\label{cor:equality}
Let $n\ge1$, and write $n=2^r u$ with $u$ odd.  Then
\begin{equation}\label{eq:equalitycriterion}
  a(n)=s(n)
  \quad\Longleftrightarrow\quad
  \text{the binary expansion of $u$ has at most one zero}.
\end{equation}
Equivalently, equality holds precisely when
\[
  u=2^m-1\quad(m\ge1),
  \qquad\text{or}\qquad
  u=2^m-1-2^j\quad(m\ge3,\ 1\le j\le m-2).
\]
\end{corollary}
\begin{proof}
By Lemma~\ref{prop:dyadic} and the Stern recurrence $s(2m)=s(m)$,
multiplication of $n$ by a power of $2$ preserves the equality $a(n)=s(n)$.
It is therefore enough to consider odd $n$.  The case $n=1$ is immediate;
assume $n>1$, and put
\[
 h:=\floor{\log_2 n},\qquad N:=n-1,\qquad w:=s_2(N).
\]
We will show that equality holds if and only if $w\ge h-1$.

Suppose first that $w\ge h-1$, and let $e$ be any hyperbinary expansion of
$N$, with total digit sum $W:=\sum_j e_j$.  Ordinary binary carrying
decreases the number of summands, so $W\ge w$.  Set $k:=W+1\ge h$.
Since $N<2^{h+1}$, the digit at position $h$ is at most $1$, and all higher
digits vanish.  Thus all digits at and above the cut $k$ are binary, and
$L_k(e)+R_k(e)=W$.  Consequently
\[
 F_e(k)=1+W-k=0.
\]
Every expansion therefore has an admissible cut.  By
Lemma~\ref{prop:cut}, the injection $\Phi_n$ is a bijection, so
$a(n)=s(n)$.

Conversely, suppose $w\le h-2$, and start with the ordinary binary
expansion $e$ of $N$.  If it has no admissible cut, then
Lemma~\ref{prop:cut} already gives $a(n)<s(n)$.  Otherwise let $k$ be
its unique admissible cut.  Since
\[
 F_e(h)=1+w-h<0,
\]
the strict decrease in \eqref{eq:Fdrop} gives $k<h$.

Let $2^j$ be the least significant nonzero term of the binary expansion of
$N$.  Because $N$ is even, $j\ge1$.  Replacing this term by two copies of
$2^{j-1}$ gives a valid hyperbinary expansion $e'$: the digit at position
$j-1$ was zero.  Its only digit $2$ is at that position, so its eligible
cuts have positions $t\ge j$.  At $t=j$, the split increases $L_t$ by $2$
and decreases $R_t$ by $1$; at every $t>j$, it increases $L_t$ by $1$ and
leaves $R_t$ unchanged.  Hence
\[
 F_{e'}(t)=F_e(t)+1\qquad(t\ge j).
\]
As $n>1$ is odd, $N\ge2^h$.  In particular, $k<h$ implies
$R_{k+1}(e)\ge1$, and \eqref{eq:Fdrop} gives
\[
 F_e(k+1)=-R_{k+1}(e)-1\le-2.
\]
Since $F_e(k)=0$ and $F_e$ is strictly decreasing, its values at cuts
$t\le k$ are nonnegative, while its values at cuts $t\ge k+1$ are at most
$-2$.  It never takes the value $-1$.  Thus $F_{e'}$ cannot vanish at an
eligible cut, and $e'$ lies outside the image of $\Phi_n$.  This proves
$a(n)<s(n)$.

Finally, subtracting one from odd $n$ changes only its final binary digit,
so $w=s_2(n)-1$.  The condition $w\ge h-1$ is therefore equivalent to
$s_2(n)\ge h$, which says that the $h+1$ binary digits of $n$ contain at
most one zero.  An odd integer whose binary expansion has no zero is $2^m-1$; if the expansion
has one zero, that zero is in an interior position, giving
$2^m-1-2^j$ with $m\ge3$ and $1\le j\le m-2$.  Restoring the removed power
of $2$ gives the stated classification.
\end{proof}
\Needspace{11\baselineskip}
The integers attaining equality are sparse.

\begin{corollary}\label{cor:sharpcount}
For every $H\ge1$, exactly
\[
 \frac{H^2-H+2}{2}
\]
of the $2^{H-1}$ integers $n$ in the interval
$2^{H-1}\le n<2^H$ satisfy $a(n)=s(n)$.  Consequently, among
$1\le n<2^H$ there are exactly
\[
 \frac{H(H^2+5)}{6}
\]
such integers.  In particular, the set of $n$ for which Stern's bound is sharp
has natural density zero: its proportion among $1\le n<X$ tends to zero as $X\to\infty$.
\end{corollary}

\begin{proof}
For $H=1$ the claim is immediate.  Assume $H\ge2$.  Write $n=2^r u$ with $u$ odd, and let $m$ be the number of binary digits of
$u$.  As $n$ runs through $2^{H-1}\le n<2^H$, the possible values of $m$ are
$1,2,\ldots,H$, each occurring for exactly one value of $r=H-m$.
By Corollary \ref{cor:equality}, for $m=1$ and $m=2$ there is one possible
odd part, while for $m\ge3$ there are $m-1$ possibilities: the expansion containing only ones,
or an expansion with its unique zero in one of the $m-2$ interior positions.  Hence
the count in this interval is
\[
 1+1+\sum_{m=3}^{H}(m-1)=\frac{H^2-H+2}{2}.
\]
Summing this expression over the first $H$ such intervals gives
$H(H^2+5)/6$.  For $2^{H-1}\le X<2^H$, the proportion below $X$ is
at most $H(H^2+5)/(6\cdot2^{H-1})$, which tends to zero.
\end{proof}

The criterion recovers the Mersenne formula in A349552 and gives a nearby
infinite family.

\begin{corollary}\label{cor:families}
For $m\ge1$,
\begin{equation}\label{eq:mersenne}
 a(2^m-1)=m.
\end{equation}
More generally, for $m\ge3$ and $1\le j\le m-2$,
\begin{equation}\label{eq:onezero}
 a(2^m-1-2^j)=(m-j)(j+1)-1.
\end{equation}
The same values persist after multiplication of the arguments by arbitrary
powers of $2$.
\end{corollary}

\begin{proof}
The relevant odd integers have no zero, or exactly one zero, in their binary
expansion, so Corollary \ref{cor:equality} reduces $a(n)$ to the Stern number
$s(n)$.  For $2^m-1=(11\cdots1)_2$, the Stern recurrence gives
$s(2^m-1)=m$.

For $n=2^m-1-2^j$, its binary expansion is
\[
 \underbrace{11\cdots1}_{m-j-1}\,0\,
 \underbrace{11\cdots1}_{j}.
\]
The pair $(s(r),s(r+1))$ transforms under appending a binary
$0$ or $1$ as
\[
 (x,y)\xrightarrow{0}(x,x+y),
 \qquad
 (x,y)\xrightarrow{1}(x+y,y).
\]
After reading the first leading $1$, the pair is $(1,1)$.  Reading the
remaining $L-1$ ones, where $L:=m-j-1$, gives $(L,1)$; after the zero it is
$(L,L+1)$; and after the final $j$ ones its first coordinate is
\[
 L+j(L+1)=(m-j)(j+1)-1.
\]
Invariance under multiplication by powers of $2$ follows from Lemma \ref{prop:dyadic}.
\end{proof}

As another immediate consequence of Theorem \ref{thm:deficit}, any
upper bound for Stern's sequence transfers to $a(n)$.  For example, the
bounds for the maxima of Stern's sequence between consecutive powers of $2$
\cite{Northshield} give
\[
 a(n)=O\!\left(n^{\log_2\varphi}\right),
 \qquad \varphi:=\frac{1+\sqrt5}{2}.
\]
A more precise description of the maximum of $a(n)$ on these intervals
is one of the questions discussed in Section~\ref{sec:discussion}.

\section{A recurrence and stabilization}\label{sec:divide}

To count recursively, we remove the first part and count the possible
remaining chains.  When the starting integer is odd, the required sum of
the remaining parts need not equal their starting integer.  We therefore
introduce counts in which the required total is allowed to vary.

\subsection{Allowing the total to vary}
For $n\ge1$ and $r\in\mathbb Z$, let $A_r(n)$ be the number of finite
\emph{halving chains} $P=(p_1,\ldots,p_k)$ of positive integers satisfying
\[
 p_1\in\mathcal H(n),\qquad
 p_i\in\mathcal H(p_{i-1})\quad(2\le i\le k),
\]
whose total is
\[
 |P|:=p_1+\cdots+p_k=n+r.
\]
The empty chain is allowed and has total $0$.  Thus
\begin{equation}\label{eq:Acentral}
 A_0(n)=a(n).
\end{equation}
The integer $r$ is the difference between the chain total and $n$.

\begin{theorem}\label{thm:Arec}
For $n\ge1$,
\begin{equation}\label{eq:A-boundary}
 A_r(n)=0\quad(r<-n),\qquad A_{-n}(n)=1,
\end{equation}
and
\begin{equation}\label{eq:A-one}
 A_r(1)=1\qquad(r\ge-1).
\end{equation}
For $m\ge1$, whenever the requested total is positive,
\begin{align}
 A_r(2m)&=A_r(m) &&(r>-2m),\label{eq:A-even}\\
 A_r(2m+1)&=A_{r+1}(m)+A_{r-1}(m+1)
 &&(r>-(2m+1)).\label{eq:A-odd}
\end{align}
In particular, these recurrences compute $a(n)=A_0(n)$.
\end{theorem}

\begin{proof}
The boundary conditions record impossibility of a negative total and the
unique empty chain of total $0$.  Starting from $1$, every nonempty chain
consists only of $1$'s, which proves \eqref{eq:A-one}.

For $2m$, every nonempty chain begins with the unique first part $m$.  Removing
that part leaves a chain from $m$ whose total must be
$(2m+r)-m=m+r$, proving \eqref{eq:A-even}.  For $2m+1$, the first part is
$m$ or $m+1$.  In the first case the remaining total is $m+1+r$, counted by
$A_{r+1}(m)$; in the second it is $m+r$, counted by $A_{r-1}(m+1)$.
This proves \eqref{eq:A-odd}.
\end{proof}

Store each computed value $A_r(m)$ so that repeated requests do not
repeat the calculation.  Starting from $A_0(n)$, after $j$ successive
recursive steps the argument $m$ is one of
$\floor{n/2^j}$ and $\ceil{n/2^j}$, while $|r|\le j$.
There are $O(\log n)$ possible depths, so at most $O((\log n)^2)$
different pairs $(m,r)$ are needed.  This counts the required values;
it does not include the cost of arithmetic on integers of increasing size.
The index $r$ also appears directly in the balance identity.  If a
nonempty chain of length $k$ is counted by $A_r(n)$ and has digits $q_i$ as in
\eqref{eq:qdef}, the proof of Lemma~\ref{lem:balance} gives the generalized
identity
\begin{equation}\label{eq:balance-r}
 \sum_{i=1}^k q_i=k-p_k-r.
\end{equation}
Encode the chain as a hyperbinary expansion $e$ of $n-1$ using
\eqref{eq:qbinary} and the ordinary binary expansion of $p_k-1$.
Then \eqref{eq:balance-r} gives $F_e(k)=-r$; the original halving
partitions are precisely the chains counted when $r=0$.

\subsection{Stabilization at Stern's sequence}
Induction on $n$ in Theorem~\ref{thm:Arec}, for all $r$ at once, gives
$0\le A_r(n)\le s(n)$: the boundary values satisfy the bound, and each
nonempty-chain recurrence is bounded by the corresponding Stern recurrence.
For a fixed starting integer, this upper bound is attained at every
sufficiently large total.  Count chains by their total using the generating function
\begin{equation}\label{eq:Bdef}
 B_n(u):=\sum_{r\ge-n}A_r(n)u^{n+r}\qquad(n\ge1).
\end{equation}
Thus $A_r(n)=[u^{n+r}]B_n(u)$, with negative-degree coefficients taken
as zero.  Set $B_0(u):=1$ separately.
\Needspace{8\baselineskip}
The empty chain contributes $1$, and removing the first part of a nonempty
chain gives
\begin{equation}\label{eq:Brec}
\begin{aligned}
 B_1(u)&=\frac1{1-u},\\
 B_{2m}(u)&=1+u^mB_m(u),\\
 B_{2m+1}(u)&=1+u^mB_m(u)+u^{m+1}B_{m+1}(u)
 \qquad(m\ge1).
\end{aligned}
\end{equation}

\Needspace{15\baselineskip}
\begin{theorem}\label{thm:stabilization}
Put $U(1):=0$, and for $n\ge2$ put
\begin{equation}\label{eq:chain-threshold}
 U(n):=n+\ceil{\log_2n}-2-s_2(n-1).
\end{equation}
Then $P_n(u):=(1-u)B_n(u)$ is a monic polynomial (its leading
coefficient is $1$) of degree $U(n)$, and $P_n(1)=s(n)$.  Consequently,
\begin{equation}\label{eq:chain-stabilization}
 A_r(n)=s(n)\quad\text{if }n+r\ge U(n),
 \qquad A_{U(n)-n-1}(n)=s(n)-1.
\end{equation}
Thus $U(n)$ is the least total from which the chain counts remain constant.
\end{theorem}

\begin{proof}
The initial series give $P_1(u)=P_2(u)=1$.  Multiplying
\eqref{eq:Brec} by $1-u$ gives
\begin{equation}\label{eq:P-rec}
\begin{aligned}
 P_{2m}(u)&=1-u+u^mP_m(u),\\
 P_{2m+1}(u)&=1-u+u^mP_m(u)+u^{m+1}P_{m+1}(u)
 \qquad(m\ge1).
\end{aligned}
\end{equation}
For $n\ge2$, successive ceiling halvings reach $1$ after
$\ell:=\ceil{\log_2n}$ steps.  Their total before that first $1$ is
\begin{align*}
 \sum_{i=1}^{\ell-1}\ceil{\frac n{2^i}}
 &=\sum_{i=1}^{\ell-1}
   \left(\floor{\frac{n-1}{2^i}}+1\right)\\
 &=n+\ell-2-s_2(n-1)=U(n),
\end{align*}
by the binary floor-sum identity \eqref{eq:floor-sum}.
This sum is nondecreasing in $n$: while $\ell$ is fixed each summand is
nondecreasing, and when $\ell$ increases an additional positive summand
appears.  Following the first halving step gives
\[
 U(2m)=m+U(m)\quad(m\ge2),\qquad
 U(2m+1)=m+1+U(m+1)\quad(m\ge1).
\]
Induct on $n$.  In the even recurrence with $m\ge2$, the shifted polynomial
has degree $m+U(m)\ge2$, so $1-u$ does not affect its leading term.
In the odd recurrence, monotonicity of $U$ gives
$m+1+U(m+1)>m+U(m)$, and the larger degree is again at least $2$.
Thus $P_n$ is monic of degree $U(n)$.
Evaluating \eqref{eq:P-rec} at $u=1$ gives Stern's recurrence, proving
$P_n(1)=s(n)$.

Write $P_n(u)=\sum_{j=0}^{U(n)}c_j u^j$.  Since $B_n=P_n/(1-u)$,
its coefficient of $u^T$ is $\sum_{j=0}^{T}c_j$, with $c_j=0$ beyond
the degree.  This equals $s(n)$ for $T\ge U(n)$, and equals $s(n)-1$
at $T=U(n)-1$, since $c_{U(n)}=1$.
When $U(n)=0$, the latter coefficient has negative degree and is zero,
as required.
\end{proof}

The threshold concerns equality at every later total; the count can
already equal $s(n)$ at some smaller totals.  For example,
$B_4(u)=1+u^2/(1-u)$ has coefficients $1,0,1,1,\ldots$.
The threshold has a particularly simple binary interpretation.

\Needspace{12\baselineskip}
\begin{corollary}\label{cor:binary-stabilization}
Let $n=2^v q$, where $q>1$ is odd, and let $z(q)$ be the number of zero
digits in the binary expansion of $q$.  Then
\begin{equation}\label{eq:binary-stabilization}
 A_r(n)=s(n)\quad(r\ge z(q)-1),
 \qquad A_{z(q)-2}(n)=s(n)-1.
\end{equation}
In particular, if the odd part of $n$ has exactly two binary zeros, then
$d(n)=1$.
\end{corollary}

\begin{proof}
Since $q>1$ is odd, its binary length is $\ceil{\log_2q}$ and
$s_2(q-1)=s_2(q)-1$.  Also
$s_2(n-1)=v+s_2(q-1)$.  Substituting these identities into
\eqref{eq:chain-threshold} gives $U(n)-n=z(q)-1$.
Now apply Theorem~\ref{thm:stabilization}, taking $r=0$ when $z(q)=2$.
\end{proof}

For example, $9=(1001)_2$ has two zero digits, so $a(9)=s(9)-1=3$,
whereas $A_r(9)=4$ for every $r\ge1$.
This extends the equality cases of Corollary~\ref{cor:equality}
to odd parts with one more zero in their binary expansion.

The formulas \eqref{eq:Drdef}--\eqref{eq:d-rec} in
Appendix~\ref{sec:appendix} give recurrences for the differences
$s(n)-A_r(n)$.  The same appendix collects the functions $B_n(u)$ in a
two-variable series.  Selecting terms with equal exponents in the two
variables recovers the generating function of $a(n)$; see
\eqref{eq:calBfunc}--\eqref{eq:a-diagonal}.

\Needspace{6\baselineskip}
\section{Distribution by length}\label{sec:fences}
\subsection{Lengths modulo three}
The counts of hyperbinary expansions in the three length classes modulo
three are as equal as possible.  We separate these classes by evaluating
$H_m(z)$ at the cube roots of unity.  This method of selecting congruence
classes of coefficients is called \emph{series multisection}
\cite{Riordan}.  The evaluation in \eqref{eq:omega-unit} follows from
Kogiso, Miyamoto, Ren, Wakui, and Yanagawa
\cite[Theorem~7.1 and Corollary~7.2]{KogisoEtAl}, together with the
hyperbinary interpretation in
\cite[Theorem~2.3]{McConvilleProppSagan}.  We give a direct proof from
\eqref{eq:Hrec}.  Let $\omega:=(-1+\sqrt{-3})/2$, so that
$\omega^3=1$ and $1+\omega+\omega^2=0$.  For $r\in\{0,1,2\}$ put
\begin{equation}\label{eq:trisection-def}
 N_r(m):=\sum_{\substack{j\ge0\\j\equiv r\pmod 3}}
          \coef{j}H_m(z).
\end{equation}

\begin{theorem}\label{prop:trisection}
For every $m\ge0$,
\begin{equation}\label{eq:omega-unit}
 H_m(\omega)\in\{0,\pm1,\pm\omega,\pm\omega^2\}.
\end{equation}
Consequently the three integers $N_0(m),N_1(m),N_2(m)$ differ pairwise by at
most one.  More precisely, if $s(m+1)=3q+\epsilon$, where
$\epsilon\in\{0,1,2\}$, then their multiset is
\[
 \{N_0(m),N_1(m),N_2(m)\}=
 \begin{cases}
  \{q,q,q\},&\epsilon=0,\\
  \{q+1,q,q\},&\epsilon=1,\\
  \{q+1,q+1,q\},&\epsilon=2.
 \end{cases}
\]
In particular,
\begin{equation}\label{eq:phi3-div}
 (1+z+z^2)\mid H_m(z)\quad\Longleftrightarrow\quad 3\mid s(m+1).
\end{equation}
\end{theorem}

\begin{proof}
Put $U_m:=H_m(\omega)$ and $V_m:=(U_m,U_{m+1})$.  If $V_m=(x,y)$,
then \eqref{eq:Hrec} and $\omega^3=1$ give
\begin{equation}\label{eq:omega-pair}
 V_{2m+1}=(\omega x,y+\omega^2x),\qquad
 V_{2m+2}=(y+\omega^2x,\omega y).
\end{equation}
Every coordinate belongs to
$\mathbb Z[\omega]:=\{a+b\omega:a,b\in\mathbb Z\}$.
For such a number $u$, write $\mathcal N(u):=u\overline u=|u|^2$,
where the bar denotes complex conjugation, and set
\[
 \mathcal Q(x,y):=\mathcal N(x)+\mathcal N(y)
  +\mathcal N(x+\omega y)+\mathcal N(x-\omega^2y).
\]
Each transformation in \eqref{eq:omega-pair} preserves $\mathcal Q$.
Indeed, for the first transformation the four arguments inside the norms become,
up to multiplication by numbers of absolute value $1$,
\[
 x,\qquad x+\omega y,\qquad x-\omega^2y,\qquad y,
\]
and the second transformation gives the same four arguments in another order.
Since $V_0=(1,\omega)$ and $\mathcal Q(1,\omega)=3$, every $V_m=(x,y)$
satisfies $\mathcal Q(x,y)=3$.  If $\mathcal N(x)=3$, all the other norms
would have to vanish, forcing $y=0$ and then
$\mathcal N(x+\omega y)=3$, a contradiction.  Nor can $\mathcal N(a+b\omega)$ equal $2$, since
\[
 \mathcal N(a+b\omega)=a^2-ab+b^2\equiv(a+b)^2\pmod3.
\]
Hence $\mathcal N(x)\in\{0,1\}$, and similarly for $y$.
The integer solutions of $a^2-ab+b^2=1$ give precisely the six numbers
$\pm1,\pm\omega,\pm\omega^2$, proving \eqref{eq:omega-unit}.

To select exponents congruent to $r$ modulo $3$, note that
$1+\omega^{j-r}+\omega^{2(j-r)}$ is $3$ when $j\equiv r\pmod3$
and $0$ otherwise.  Applying this identity to each coefficient gives
\begin{equation}\label{eq:trisection-filter}
 N_r(m)=\frac13\left(H_m(1)+\omega^{-r}H_m(\omega)
              +\omega^{-2r}H_m(\omega^2)\right).
\end{equation}
Since $H_m$ has integer coefficients,
$H_m(\omega^2)=\overline{H_m(\omega)}$, while directly from
\eqref{eq:trisection-def},
\[
 H_m(\omega)=(N_0-N_2)+(N_1-N_2)\omega.
\]
The six possible nonzero values therefore give, according to
$H_m(1)=s(m+1)$, the permutations of $(q+1,q,q)$ or $(q+1,q+1,q)$ displayed
above; zero gives $(q,q,q)$.  Finally, $H_m(\omega)=0$ is equivalent to
$1+z+z^2$ dividing $H_m(z)$, and by the same classification this occurs exactly
when $3\mid H_m(1)=s(m+1)$.
\end{proof}

\subsection{Fence posets}

A description by partially ordered sets explains the distribution of
lengths.  Earlier refinements of Stern polynomials count hyperbinary
expansions by various properties of their digits
\cite{KlavzarMilutinovicPetr,DilcherEricksen,MansourShattuck}.
Brunetti and D'Aniello arranged the vertices of the graph of hyperbinary
expansions by total digit sum \cite{BrunettiDaniello}.
McConville, Propp, and Sagan study the same length polynomial using
partially ordered sets, or \emph{posets}.  Order hyperbinary partitions by
\emph{refinement}, with finer partitions below coarser ones.  An elementary
refinement replaces $2^j$ by two copies of $2^{j-1}$ whenever the result
is still hyperbinary.  Their result identifies this ordering with
inclusion among order ideals of a fence \cite{McConvilleProppSagan}.

Here are the definitions needed for that description.  An element $y$
\emph{covers} $x$ if $x<y$ and no element lies strictly between them.
The \emph{Hasse diagram} joins each such pair, with the larger element
above the smaller.  A \emph{fence} is a finite poset whose Hasse diagram,
with edge directions ignored, is a path; the edge directions may vary
along the path.  We also allow the empty fence.  An \emph{order ideal}
is a subset $I$ such that $x\in I$ and $y\le x$ imply $y\in I$.
Write $J(F)$ for the order ideals of $F$, ordered by inclusion.
Intersections and unions of order ideals are again order ideals, and
these operations distribute over each other.  Thus $J(F)$ is a
\emph{distributive lattice}.

For $m\ge1$, write its binary expansion as $b_1\cdots b_t$, with $b_1=1$.
Its \emph{principal prefix} $b_1\cdots b_r$ consists of the digits before
the rightmost zero; if all digits are one, the prefix is empty.
Following McConville, Propp, and Sagan, let $\mathcal F(m)$ have elements
$x_1,\ldots,x_r$, with covering relations only between adjacent elements:
\[
 x_i<x_{i-1}\ \text{if }b_i=0,\qquad
 x_i>x_{i-1}\ \text{if }b_i=1\qquad(2\le i\le r).
\]
Their Theorems~3.16--3.17 identify the hyperbinary refinement poset with
$J(\mathcal F(m))$ and relate length to the size of an ideal
\cite{McConvilleProppSagan}.  If $r:=|\mathcal F(m)|$, their length formula may be written in our notation as
\begin{equation}\label{eq:fence-H}
 H_m(z)=z^{s_2(m)} C_m(z),\qquad
 C_m(z):=\sum_{I\in J(\mathcal F(m))}z^{r-|I|}.
\end{equation}
Thus the length of an expansion is $s_2(m)$ plus the number of fence
elements omitted from the corresponding ideal.  The latter number,
$r-|I|$, is often called the \emph{corank} of $I$ in $J(\mathcal F(m))$.
Theorem~\ref{prop:trisection} therefore gives a result about all fences,
without reference to hyperbinary expansions.

\begin{corollary}\label{cor:fence-trisection}
For a finite fence $F$, put
\[
 \nu_j(F):=\#\{I\in J(F):|I|\equiv j\pmod3\},\qquad j=0,1,2.
\]
Then $\nu_0(F),\nu_1(F),\nu_2(F)$ differ pairwise by at most one.
\end{corollary}

\begin{proof}
The empty fence is immediate.  For a fence with $r\ge1$ elements, record its
successive edge directions by digits $1b_2\cdots b_r$, using the
rule for covering relations just defined.  Choose $m$ whose
binary expansion is $1b_2\cdots b_r01$.  Its principal prefix is
$1b_2\cdots b_r$, so the construction of McConville, Propp, and Sagan gives
$\mathcal F(m)\cong F$.  By \eqref{eq:fence-H}, multiplication by
$z^{s_2(m)}$ only cyclically shifts the three coefficient classes modulo $3$,
and the exponent $r-|I|$ in $C_m(z)$ only reverses and shifts the residue class
of $|I|$.  The result now follows from Theorem~\ref{prop:trisection}.
\end{proof}

For $m=12=(1100)_2$, for example, $\mathcal F(12)$ has three elements:
two minimal elements covered by a common maximum.  Its five order ideals omit
$0,1,2,2,3$ elements, respectively, so
\[
 C_{12}(z)=1+z+2z^2+z^3,
 \qquad
 H_{12}(z)=z^2C_{12}(z).
\]
The coefficients $1,1,2,1$ in Example~\ref{ex:13} therefore count
these ideals by the number of elements they omit.

Combining \eqref{eq:longlengths} with \eqref{eq:fence-H} gives a structural
interpretation of every sufficiently long halving partition.

\begin{corollary}\label{prop:fence}
Let $n\ge2$, put $h:=\floor{\log_2n}$, $m:=n-1$, and $r:=|\mathcal F(m)|$.  For every $k\ge h$,
\begin{equation}\label{eq:fence-ak}
 a_k(n)=
 \#\left\{I\in J(\mathcal F(m)):
 r-|I|=k-1-s_2(m)\right\}.
\end{equation}
\end{corollary}

A finite sequence is \emph{unimodal} if it weakly increases up to some
position and weakly decreases thereafter.  Kantarc{\i} O\u{g}uz and
Ravichandran proved that, for every fence $F$, the numbers of its order
ideals of sizes $0,1,\ldots,|F|$ form a unimodal sequence
\cite{OguzRavichandran}.  These are the coefficients of the
\emph{rank polynomial} $\sum_{I\in J(F)}z^{|I|}$.
Reversing a unimodal sequence preserves unimodality, so
\eqref{eq:fence-H} shows that the coefficients of every $H_m(z)$
are unimodal.  The counts for sufficiently long halving partitions form
a consecutive portion of that sequence.

\begin{corollary}\label{cor:unimodal}
For $n\ge2$, let $h:=\floor{\log_2n}$.  Then
\[
 a_h(n),a_{h+1}(n),\ldots,a_{L(n)}(n)
\]
is unimodal.
\end{corollary}

\begin{proof}
By \eqref{eq:longlengths}, this sequence is the consecutive coefficient segment
\[
 \coef{h-1}H_{n-1}(z),\ \coef{h}H_{n-1}(z),\ldots,
 \coef{L(n)-1}H_{n-1}(z).
\]
The coefficient sequence of $H_{n-1}(z)$ is unimodal by
\eqref{eq:fence-H} and the theorem of Kantarc{\i} O\u{g}uz and Ravichandran.  A consecutive
segment of a unimodal sequence is unimodal.
\end{proof}

\Needspace{16\baselineskip}
\begin{example}\label{ex:nonunimodal}
Take $n:=49185=3\cdot2^{14}+2^5+1$.  Then $h=15$, $L(n)=19$, and
\begin{equation}\label{eq:nonunimodal-profile}
 (a_{13}(n),a_{14}(n),\ldots,a_{19}(n))=(6,5,9,7,5,3,1).
\end{equation}
Here is an exact check using only \eqref{eq:Hrec}.  Put
$S_j:=z+z^2+\cdots+z^j$.  Repeated application of that recurrence gives
\[
\begin{aligned}
 H_{8224}(z)&=zS_8+zS_5S_9,\\
 H_{16416}(z)&=zS_9+zS_5S_{10},\\
 H_{49184}(z)&=z^3\bigl(1+S_5+(1+z)S_8+(1+z)S_5S_9\bigr).
\end{aligned}
\]
Theorem~\ref{thm:main} now gives
$a_{13}(n)=\coef{7}H_{8224}(z)=6$ and
$a_{14}(n)=\coef{11}H_{16416}(z)=5$.
For $k\ge15$, use the coefficients of $H_{49184}$ in degrees $k-1$.
In degrees $14$ through $18$, only the last term in the displayed
expression contributes, giving $9,7,5,3,1$.
The upper bound on $n$ in \eqref{eq:fixedlength-support} excludes $k\le12$,
since $11\cdot2^{12}+1=45057<n$, and Theorem~\ref{thm:longest} excludes
$k>19$.  Thus \eqref{eq:nonunimodal-profile} lists all the nonzero counts by length.  It decreases from $6$ to $5$ and then increases to $9$;
the counts for $k\ge h$ are nevertheless decreasing.
\end{example}

\section{Discussion}\label{sec:discussion}
The balance identity connects a local rounding rule with Stern's sequence
and order ideals of fences.  The cut criterion identifies
the hyperbinary expansions that satisfy the required total.
Theorem~\ref{thm:stabilization} gives a complementary role for Stern's
sequence as the eventual chain count when the prescribed total increases.
The binary equality criterion and Corollary~\ref{cor:binary-stabilization}
determine the deficit when the odd part has at most two zero digits.
For odd parts with more zero digits, is there a direct formula for the
deficit in terms of those digits?  More generally, how large can $a(n)$
be among integers with a prescribed number of binary digits?

The recurrence in Section~\ref{sec:divide} requires only
$O((\log n)^2)$ values $A_r(m)$ for each input $n$, but the range of $r$
is not bounded independently of $n$.  This raises the question of
whether $a(n)$ is \emph{$2$-regular}: does there exist a finite family
of integer sequences containing $a(n)$ such that, for each sequence
$b(n)$ in the family, both $b(2n)$ and $b(2n+1)$ are integer linear
combinations of sequences in the family \cite{HeubergerKrennLechner}?

Weighting occurrences of the digit $1$ in unrestricted hyperbinary
expansions by $t$ gives the classical Stern polynomials
\cite[Theorem~2]{KlavzarMilutinovicPetr}.  A refinement specific to this
problem would retain \eqref{eq:balance} while assigning weight $t$ to
each occurrence of $q_i=1$, and ask how the sum condition changes the
resulting polynomials.  One may also replace halving by division by a
larger base.  Related hyper-$m$-ary expansions use powers of that base
with bounded multiplicities \cite{MansourShattuckMary}.  Which of the
correspondences and exact bounds extend to this setting?

\section*{Acknowledgments}
ChatGPT was used to assist with developing and checking arguments and
calculations, and with revising the exposition.  The author takes
responsibility for the content of the paper.

\Needspace{12\baselineskip}
\begingroup
\small
\raggedright

\endgroup

\clearpage
\appendix
\numberwithin{equation}{section}
\renewcommand{\theequation}{\thesection\arabic{equation}}
\makeatletter
\renewcommand{\@seccntformat}[1]{Appendix~\csname the#1\endcsname.\hskip .5em}
\makeatother

\section{Supplementary formulas and a comparison}\label{sec:appendix}
These formulas supplement the fixed-length count and the recurrences of
Section~\ref{sec:divide}.  They give an alternative coefficient formula
for short lengths, recurrences for the differences $s(n)-A_r(n)$, and
a two-variable generating function for chains.  We then give the detailed
comparison with the payphone interval calculation introduced in
Section~\ref{sec:intro}.

\subsection{Alternative formulas for short lengths}\label{app:short-length}

For $1\le k<h:=\floor{\log_2 n}$, put $r_k:=n\bmod 2^k$.  Since
$M_k(n)=2^k+r_k-1$ and $Q_k(n)=(n-r_k)/2^k$, the fixed-length
formula~\eqref{eq:main} gives
\begin{align}
 a_k(n)
 &=[x^{\,2^k+r_k-1}z^{\,k2^k+r_k-n}]G(x,z^{2^k}) \notag\\
 &=[x^{\,2^k+r_k-1}z^{\,(k-1)2^k-n+1}]G(x/z,z^{2^k}).
 \label{eq:rescaled-ak}
\end{align}
The substitution is interpreted separately for each coefficient of $x$;
that coefficient may contain negative as well as nonnegative powers of $z$.  The second form confines the residue $r_k$ to the
exponent of $x$.

For $n\ge2$, substituting \eqref{eq:main} into
\eqref{eq:deficit} also gives, with $J(h)$ defined in \eqref{eq:Jdef},
\begin{equation}\label{eq:deficit-explicit}
 d(n)=\coef{h-2}\frac{H_{n-1}(z)}{1-z}
 -\sum_{j=1}^{J(h)}\coef{h-j-Q_{h-j}(n)}H_{M_{h-j}(n)}(z).
\end{equation}

\subsection{Recurrences for the deficit}\label{app:deficit-states}
For $n\ge1$ and $r\in\mathbb Z$, define the difference between Stern's
number and the chain count by
\begin{equation}\label{eq:Drdef}
 D_r(n):=s(n)-A_r(n).
\end{equation}
Then $D_0(n)=d(n)$, and Stern's recurrences together with
\eqref{eq:A-even} and \eqref{eq:A-odd} give
\begin{align}
 D_r(2m)&=D_r(m),\label{eq:D-even}\\*
 D_r(2m+1)&=D_{r+1}(m)+D_{r-1}(m+1),\label{eq:D-odd}
\end{align}
for $m\ge1$, with $r>-2m$ in \eqref{eq:D-even} and
$r>-(2m+1)$ in \eqref{eq:D-odd}.  The boundary values
are explicit:
\[
 D_r(n)=s(n)\quad(r<-n),\qquad
 D_{-n}(n)=s(n)-1,
\]
and $D_r(1)=0$ for $r\ge-1$.  In particular,
\begin{equation}\label{eq:d-rec}
 d(2m)=d(m),\qquad
 d(2m+1)=D_1(m)+D_{-1}(m+1).
\end{equation}
The even case uses $d$ alone; the odd case requires the auxiliary counts
with indices $1$ and $-1$.
Theorem~\ref{thm:stabilization} supplies the further values
$D_r(n)=0$ whenever $n+r\ge U(n)$.

\Needspace{12\baselineskip}
\subsection{A bivariate series and its diagonal}\label{app:chain-sum}
Using the generating function $B_n(u)$ for chain totals from \eqref{eq:Bdef}, put
\begin{equation}\label{eq:calBdef}
 \mathcal B(x,u):=\sum_{n\ge1}B_n(u)x^n.
\end{equation}
Separating even and odd indices in \eqref{eq:Brec} gives
\begin{equation}\label{eq:calBfunc}
 \mathcal B(x,u)=\frac{x}{1-x}
 +(x^{-1}+1+x)\mathcal B(ux^2,u).
\end{equation}
The factor $x^{-1}$ introduces no negative powers of $x$, since
$\mathcal B(ux^2,u)$ is divisible by $x^2$.
The \emph{diagonal} selects coefficients whose two exponents agree:
\[
 \operatorname{Diag}_{x,u}\sum_{i,j}b_{ij}x^iu^j:=\sum_n b_{nn}t^n.
\]
Then \eqref{eq:Acentral} and \eqref{eq:Bdef} give
\begin{equation}\label{eq:a-diagonal}
 \sum_{n\ge1}a(n)t^n=\operatorname{Diag}_{x,u}\mathcal B(x,u).
\end{equation}
Thus the ordinary generating function is the diagonal of the series
specified by \eqref{eq:calBfunc}.

\Needspace{12\baselineskip}
\subsection{The payphone comparison}\label{sec:payphone}
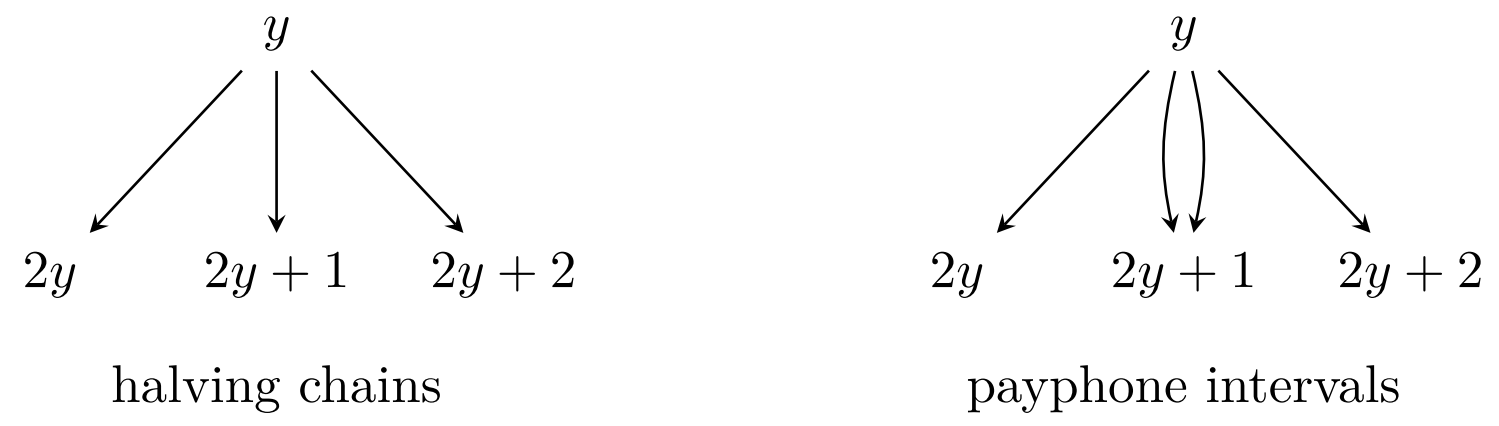
\begin{figure}[ht]
\centering
\begin{tikzpicture}[>=stealth, node distance=1.35cm and 1.6cm, line width=0.5pt, every node/.style={font=\small}]
  \node (h0) at (0,1.6) {$y$};
  \node (h1) at (-1.5,0) {$2y$};
  \node (h2) at (0,0) {$2y+1$};
  \node (h3) at (1.5,0) {$2y+2$};
  \draw[->] (h0) -- (h1);
  \draw[->] (h0) -- (h2);
  \draw[->] (h0) -- (h3);
  \node at (0,-0.75) {halving chains};

  \node (p0) at (6,1.6) {$y$};
  \node (p1) at (4.5,0) {$2y$};
  \node (p2) at (6,0) {$2y+1$};
  \node (p3) at (7.5,0) {$2y+2$};
  \draw[->] (p0) -- (p1);
  \draw[->,bend left=13] (p0) to (p2);
  \draw[->,bend right=13] (p0) to (p2);
  \draw[->] (p0) -- (p3);
  \node at (6,-0.75) {payphone intervals};
\end{tikzpicture}
\caption{Local moves for halving chains (left) and payphone interval
multiplicities (right).  The middle move has multiplicity one or two,
respectively.}
\label{fig:trees}
\end{figure}

Shift the parts of a halving chain by one, $y_i:=p_i-1$.
Equation~\eqref{eq:qdef} becomes
\begin{equation}\label{eq:reverse-tree}
 y_{i-1}=2y_i+q_i,\qquad q_i\in\{0,1,2\}.
\end{equation}
Read backwards, this gives a path from $y_k$ to $y_0$ with steps
$y\mapsto2y,2y+1,2y+2$.  Regard these steps as edges of a directed graph
on the nonnegative integers.  Every step from a positive integer
strictly increases it, but different paths can end at the same integer:
for example, $1\to4$ and $1\to2\to4$.  Treating different paths from
$1$ as different vertices gives a rooted tree whose integer labels may
repeat.  In the payphone interval calculation, the middle step is counted
twice, as in Figure~\ref{fig:trees}.  The loop at $0$ records repetitions
of the terminal part $1$ and does not occur in the tree starting from $1$.

Expansions using the digits $0,1,2$ are called \emph{quasi-binary} in the
payphone paper \cite{AlekseyevPayphone}.  A path of length $\ell$ starting at $y$, with
successive edge digits $e_{\ell-1},\ldots,e_0$, ends at
\[
 2^\ell y+\sum_{j=0}^{\ell-1}e_j2^j.
\]
The possible values of $\sum_j e_j2^j$, counted with one copy of each
step, are recorded by
\begin{equation}\label{eq:hyper-product}
 \prod_{j=0}^{\ell-1}\left(1+x^{2^j}+x^{2^{j+1}}\right).
\end{equation}
Giving the middle move two realizations replaces each factor by
$(1+x^{2^j})^2$, and the product becomes
\begin{equation}\label{eq:payphone-product}
 \prod_{j=0}^{\ell-1}(1+2x^{2^j}+x^{2^{j+1}})
 =\left(\frac{1-x^{2^\ell}}{1-x}\right)^2.
\end{equation}
Its coefficients rise from $1$ to $2^\ell$ and then fall symmetrically.
The coefficients of \eqref{eq:hyper-product} stabilize, as
$\ell$ increases, to the Stern coefficients in \eqref{eq:stern-product}.

Here is the precise role of \eqref{eq:payphone-product} in the payphone
calculation.  In an interval flanked by occupied phones, selecting a
middle phone splits an interval of length $j$ into lengths
$\floor{(j-1)/2}$ and $\ceil{(j-1)/2}$.
Let $S(N)$ be the multiset of lengths of an initial interval of length $N$
and all its nonempty descendants under repeated splitting.
Set $S(0):=\varnothing$.
If $\#_yS(N)$ denotes the multiplicity of $y$, then
\cite[Lemma~1]{AlekseyevPayphone} gives, for $N,y\ge1$,
\begin{equation}\label{eq:payphone-intervals}
 \#_yS(N)=\sum_{\ell\ge0}[x^{\,N-2^\ell y}]
 \left(\frac{1-x^{2^\ell}}{1-x}\right)^2.
\end{equation}
These interval multiplicities enter the permutation formulas together
with factors counting the possible orders of interval selection
\cite[Section~4]{AlekseyevPayphone}.
Thus the local products compare ingredients in the two enumerations;
the halving count also requires the balance identity \eqref{eq:balance}.

More generally, giving the middle move weight $t$ produces the classical
Stern polynomials, denoted here by $\mathsf S_n(t)$:
\begin{equation}\label{eq:classical-stern-polynomials}
 \sum_{N\ge0}\mathsf S_{N+1}(t)x^N
 =\prod_{j\ge0}(1+t x^{2^j}+x^{2^{j+1}}).
\end{equation}
They enumerate hyperbinary expansions of $n-1$ by the number of digits
$1$, and satisfy $\mathsf S_n(1)=s(n)$ and $\mathsf S_n(2)=n$
\cite[Theorem~2]{KlavzarMilutinovicPetr}.
The choices $t=1$ and $t=2$ give the two products above, without imposing
the sum condition for halving partitions.  Retaining that condition gives
the refinement mentioned in Section~\ref{sec:discussion}.

Finally, the payphone interval count $f(n):=\#_1S(n-1)$ satisfies
$f(1)=0$, $f(2)=1$, and
\[
 f(n)=f(\floor{n/2})+f(\ceil{n/2})\qquad(n\ge3).
\]
This is a recurrence that splits its argument into its floor and ceiling
halves, of the kind studied by Hwang, Janson, and Tsai
\cite{HwangJansonTsai}.  Theorem~\ref{thm:Arec} also halves the argument,
while retaining the difference $r$ between the chain total and its
starting integer.

\end{document}